\documentclass{amsart}

\newcommand{\supp}{\mbox{supp}\,}
\renewcommand{\span}{\mbox{span}\,}
\renewcommand{\r}{{\Bbb R}}
\newcommand{\z} {{\Bbb Z}}
\newcommand{\n} {{\Bbb N}}
\newcommand{\ddd}{,\dots,}
\newcommand{\lll}{\left(}
\newcommand{\rrr}{\right)}
\newcommand{\ex}[1]{e^{2\pi i{#1}}}

\renewcommand{\phi}{{\varphi}}

\newcommand{\cD}{{\mathcal D}}

\newcommand{\bN}{{\mathbb N}}

\newcommand{\bQ}{{\mathbb Q}}
\newcommand{\bR}{{\mathbb R}}
\newcommand{\bC}{{\mathbb C}}

\newcommand{\be}{\begin{equation}}
\newcommand{\ee}{\end{equation}}
\newcommand{\ba}{\begin{eqnarray}}
\newcommand{\ea}{\end{eqnarray}}
\newcommand{\ban}{\begin{eqnarray*}}
\newcommand{\ean}{\end{eqnarray*}}

\newtheorem{theorem}{Theorem}[section]
\newtheorem{lemma}[theorem]{Lemma}

\theoremstyle{definition}
\newtheorem{definition}[theorem]{Definition}
\newtheorem{example}[theorem]{Example}

\theoremstyle{remark}
\newtheorem{remark}[theorem]{Remark}

\theoremstyle{proposition}
\newtheorem{proposition}[theorem]{Proposition}

\theoremstyle{proposition}
\newtheorem{corollary}[theorem]{Corollary}

\numberwithin{equation}{section}

\begin{document}

\title[$p$-Adic multiresolution analyses
] {$p$-Adic multiresolution analyses
}

\author{S.~Albeverio}
\address{Universit\"at Bonn, Institut f\"ur Angewandte
Mathematik, Abteilung Sto\-chas\-tik, Wegelerstra\ss e 6, D-53115 Bonn
and Interdisziplinäres Zentrum f\"ur Komplexe Systeme, Universit\"at Bonn,
R\"omerstra\ss e 164 D-53117, Bonn, Germany}
\email{albeverio@uni-bonn.de}
\thanks{The paper was supported in part by
DFG Project 436 RUS 113/951. The second author was supported in part by Grants
06-01-00471 and 07-01-00485 of RFBR.
The third author was supported in part by Grant
06-01-00457 of RFBR}

\author{S.~Evdokimov}
\address{St.-Petersburg Department of Steklov Institute  of Mathematics, St.-Petersburg,
 Fontanka-27, 191023 St. Petersburg, RUSSIA }
 \email{evdokim@pdmi.ras.ru}

\author{M.~Skopina}
\address{Department of Applied Mathematics and Control Processes,
St. Petersburg State University, \ Universitetskii pr.-35,
198504 St. Petersburg, Russia.}
\email{skopina@MS1167.spb.edu}

\subjclass[2000]{Primary  42C40, 11E95; Secondary 11F85}

\date{}


\keywords{$p$-adic multiresolution analysis; refinable equations,
wavelets.}

\begin{abstract} We study $p$-adic multiresolution analyses
(MRAs). A complete characterisation of  test   functions
generating a MRA (scaling functions) is given. We prove that only
$1$-periodic test functions may be taken as orthogonal scaling
functions and that all such scaling functions generate Haar MRA.
We also suggest a method of constructing sets of wavelet functions
and prove that any set of wavelet functions generates a  $p$-adic
wavelet frame. \end{abstract}

\maketitle

\section{Introduction}
\label{s1}

In the early nineties a general scheme for the
construction of wavelets (of real argument) was developed. This
scheme is based on the notion of multiresolution analysis (MRA in the sequel)
introduced by Y.~Meyer and S.~Mallat~\cite{Mallat-1}, \cite{Meyer-1} (see also,
e.g.,  ~\cite{31}, ~\cite{NPS}). Immediately specialists started to
implement new wavelet systems.
Nowadays it is difficult to find an engineering area where wavelets are
not applied.

 In the $p$-adic setting, the situation is as follows.
In 2002 S.~V.~Kozyrev~\cite{Koz0} found a compactly supported
$p$-adic wavelet basis for ${ L}^2(\bQ_p)$ which is an analog of
the Haar basis. It even turned out that these wavelets were
eigenfunctions of $p$-adic pseudo-differential
operators~\cite{Koz2}. J.J.~Benedetto and
R.L.~Benedetto~\cite{Ben-Ben}, \cite{Ben1}, however,
discussed if it is possible to construct other $p$-adic wavelets
with the same set of translations which are not a group. In
particular, R.L.~Benedetto~\cite[p. 28]{Ben1}
 had doubts that a MRA-theory could
 be developed because  discrete subgroups do not exist in $\bQ_p$.
Indeed, the latter seems to be an obstacle for the development of
a MRA theory.
On the other hand, A.~Khrennikov and V.~Shelkovich~\cite{Kh-Sh1}
conjectured that the equality
\begin{equation}
\label{62.0-3}
\phi(x)=\sum_{r=0}^{p-1}\phi\Big(\frac{1}{p}x-\frac{r}{p}\Big),
\quad x\in \bQ_p,
\end{equation}
may be considered as a {\it refinement equation} for the Haar MRA
generating Kozyrev's wavelets. A solution $\phi$ of this equation
({\it a refinable function}) is the characteristic function of the
unit disc. We note that equation (\ref{62.0-3}) reflects a {\it
natural} ``self-similarity'' of the space $\bQ_p$: the unit disc
$B_{0}(0)=\{x: |x|_p \le 1\}$ is represented as the union
of $p$ mutually {\it disjoint} discs
$B_{-1}(r)=\bigl\{x: |x-r|_p \le p^{-1}\bigr\}$, $r=0,\dots,p-1$.
Following this idea, the notion of $p$-adic MRA
was introduced and a general scheme for its
construction was described in~\cite{S-Sk-1}. Also, using
(\ref{62.0-3}) as a generating refinement equation, this scheme was
realized to construct the $2$-adic Haar MRA. In contrast to the real
setting, the {\it refinable function} $\phi$ generating the Haar MRA
is {\em periodic}, which implies the existence of {\em infinitly many
different} orthonormal wavelet bases in the same Haar MRA. One of
them coincides with Kozyrev's wavelet basis.
The authors of~\cite{Kh-Sh-S} described a wide class of  functions generating
a MRA, but all of these functions are $1$-periodic.
In the present paper we prove that there exist no other orthogonal
test scaling functions generating
a MRA, except for those described in~\cite{S-Sk-1}.
Also, the MRAs generated by arbitrary test scaling functions (not necessarily
orthogonal) are considered and a criterion for a test function to generate
such a MRA is found. The non-group structure of the set of standard translations
is compensated by the fact that the sample spaces are invariant with respect
to all translations by the elements of $\bQ_p$. Finally we develop a method
to construct a wavelet frame based on a given MRA.

Here and in what follows, we shall systematically use the
notation and the results from~\cite{Vl-V-Z}.
Let $\bN$, $\z$, $\r$, $\bC$ be the sets of positive integers, integers,
real numbers, complex numbers, respectively.
The field $\bQ_p$ of $p$-adic numbers is defined as the completion
of the field of rational numbers $\bQ$ with respect to the
non-Archimedean $p$-adic norm $|\cdot|_p$. This $p$-adic norm
is defined as follows: $|0|_p=0$; if $x\ne 0$, $x=p^{\gamma}\frac{m}{n}$,
where $\gamma=\gamma(x)\in \z$
and the integers $m$, $n$ are not divisible by $p$, then
$|x|_p=p^{-\gamma}$.
The norm $|\cdot|_p$ satisfies the strong triangle inequality
$|x+y|_p\le \max(|x|_p,|y|_p)$.
The canonical form of any $p$-adic number $x\ne 0$ is
\begin{equation}
\label{2}
x=p^{\gamma}\sum_{j=0}^\infty x_jp^j
\end{equation}
where $\gamma=\gamma(x)\in \z$, \ $x_j\in D_p:=\{0,1,\dots,p-1\}$, $x_0\ne 0$.
The fractional part $\{x\}_p$ of the number $x$ equals by
definition $p^{\gamma}\sum_{j=0}^{-\gamma-1} x_jp^j$.
Thus, $\{x\}_p=0$ if and only if $\gamma\ge0$. We also set $\{0\}_p=0$.

Denote by $B_{\gamma}(a)=\{x\in \bQ_p: |x-a|_p \le p^{\gamma}\}$
the disc of radius $p^{\gamma}$ with the center at a point $a\in \bQ_p$,
$\gamma \in \z$. Any two balls in $\bQ_p$ either are disjoint or one
contains the other. We observe that $B_0(0)=\{x\in \bQ_p: \{x\}_p=0\}$.

There exists the Haar measure $dx$ on $\bQ_p$ which is  positive,
invariant under the shifts, i.e., $d(x+a)=dx$, and normalized by
$\int_{|\xi|_p\le 1}\,dx=1$.
A complex-valued function $f$ defined on $\bQ_p$ is called
{\it locally-constant} if for any $x\in \bQ_p$ there exists
an integer $l(x)\in \z$ such that
$f(x+y)=f(x)$, $y\in B_{l(x)}(0)$.
Denote by ${\cD}$ the linear space of locally-constant compactly
supported functions (so-called test functions).
The space ${\cD}$ is an analog of the  Schwartz space in the real analysis.

The Fourier transform of $\varphi\in {\cD}$ is defined as
$$
{\widehat\phi}(\xi)=F[\varphi](\xi)=\int_{\bQ_p}\chi_p(\xi x)\varphi(x)\,dx,
\ \ \ \xi \in \bQ_p,
$$
where $\chi_p(\xi x)=e^{2\pi i\{\xi x\}_p}$ is the additive character for
the field $\bQ_p$, and $\{\cdot\}_p$ is the fractional part of a $p$-adic number.
The Fourier transform is a linear isomorphism taking ${\cD}$ into
${\cD}$. The Fourier transform is extended to ${ L}^2(\bQ_p)$ in a
standard way and the Plancherel equality holds
$$
\int\limits_{\bQ_p}f(x)g(x)\,dx=
\int\limits_{\bQ_p}\widehat f(\xi)\widehat g(\xi)\,d\xi,\quad f,g\in L^2(\bQ_p).
$$
If $f\in{ L}^2(\bQ_p)$, $0\ne a\in \bQ_p$, \ $b\in \bQ_p$,
then:
\begin{equation}
\label{014}
F[f(a\cdot+b)](\xi)
=|a|_p^{-1}\chi_p\Big(-\frac{b}{a}\xi\Big)F[f]\Big(\frac{\xi}{a}\Big).
\end{equation}
Besides,
\begin{equation}
\label{14.1}
F[\Omega(|\cdot|_p)](\xi)=\Omega(|\xi|_p),
\quad \xi \in \bQ_p,
\end{equation}
where $\Omega$ is the characteristic function of the interval  $[0,\,1]$.

\section{Multiresolution analysis}
\label{s2}

Let us consider the set
$$
I_p=\{a\in \bQ_p:\{a\}_p=a\}.
$$
Since $B_0(0)=\{x\in \bQ_p: \{x\}_p=0\}$,
we have the following decomposition of
$\bQ_p$ into the union of mutually  disjoint discs:
$\bQ_p=\bigcup_{a\in I_p}B_{0}(a).$
Thus, $I_p$ can be considered as a {\em ``natural'' set of translations} for $\bQ_p$.

\begin{definition}
\label{de1} \rm
A collection of closed spaces
$V_j\subset L^2(\bQ_p)$, $j\in\z$, is called a
{\it multiresolution analysis {\rm(}MRA{\rm)} in $ L^2(\bQ_p)$} if the
following axioms hold

(a) $V_j\subset V_{j+1}$ for all $j\in\z$;

(b) $\bigcup_{j\in\z}V_j$ is dense in $ L^2(\bQ_p)$;

(c) $\bigcap_{j\in\z}V_j=\{0\}$;

(d) $f(\cdot)\in V_j \Longleftrightarrow f(p^{-1}\cdot)\in V_{j+1}$
for all $j\in\z$;

(e) there exists a function $\phi \in V_0$
such that $V_0:=\overline{\span\{\phi(x-a):\ a\in I_p\}}$.
\end{definition}

The function $\phi$ from axiom (e) is called {\em scaling}.
One also says that a MRA is generated by its scaling function $\phi$
(or $\phi$ generates the MRA).
It follows  immediately from axioms (d) and (e) that
\be
V_j:=\overline{\span\{\phi(p^{-j}x-a):\ a\in I_p\}},\quad j\in \z.
\label{17}
\ee

An important class of MRAs consists of those
generated by so-called {\em orthogonal scaling functions}.
A scaling function  $\phi$ is said to be orthogonal if
 $\{\phi(x-a), a\in I_p\}$ is an orthonormal
basis for $V_0$. Consider such a MRA.
Evidently, the functions $p^{j/2}\phi(p^{-j}x-a)$,
$a\in I_p$, form an orthonormal basis for $V_j$, $j\in\z$.
According to the standard scheme (see, e.g.,~\cite[\S 1.3]{NPS})
for the construction of MRA-based wavelets, for each $j$, we define
a space $W_j$ ({\em wavelet space}) as the orthogonal complement
of $V_j$ in $V_{j+1}$, i.e., $V_{j+1}=V_j\oplus W_j$, $j\in \z$,
where $W_j\perp V_j$, $j\in \z$. It is not difficult to see that
\begin{equation}
\label{61.0}
f(\cdot)\in W_j \Longleftrightarrow f(p^{-1}\cdot)\in W_{j+1},
\quad\text{for all}\quad j\in \z
\end{equation}
and $W_j\perp W_k$, $j\ne k$.
Taking into account axioms (b) and (c), we obtain
\begin{equation}
\label{61.1}
{\bigoplus\limits_{j\in\z}W_j}= L^2(\bQ_p)
\quad \text{(orthogonal direct sum)}.
\end{equation}
If we now find   functions $\psi^{(\nu)} \in W_0$, $\nu\in A$,
such that the functions $\psi^{(\nu)}(x-a)$, $a\in I_p, \nu\in A$, form an orthonormal
basis for $W_0$, then, due to~(\ref{61.0}) and (\ref{61.1}),
the system $\{p^{j/2}\psi^{(\nu)}(p^{-j}x-a), a\in I_p, j\in\z , \nu\in A\}$
is an orthonormal basis for $ L^2(\bQ_p)$.
Such a function $\psi^{(\nu)}$ are called a {\em wavelet function} and
the basis is a {\em wavelet basis}.

Another interesting class of scaling functions consists of
functions $\phi$  for which  $\{\phi(x-a), a\in I_p\}$ is a Riesz system.
Probably, adopting the ideas developed for the real setting, one can
use MRAs generated by such functions  for constructing
dual biorthogonal wavelet systems. This topic is, however, out of our consideration
in the present paper.

In  Section~\ref{s3} we will discuss how to construct a $p$-adic wavelet frame
based on an arbitrary MRA generated by a test function.

Let $\phi$ be an orthogonal scaling function for a MRA $\{V_j\}_{j\in\z}$. Since
the system $\{p^{1/2}\phi(p^{-1}x-a), a\in I_p\}$
is a basis for $V_1$ in this case, it follows from axiom (a) that
\begin{equation}
\label{62.0-2*}
\phi(x)=\sum_{a\in I_p}\alpha_a\phi(p^{-1}x-a),
\quad \alpha_a\in \bC.
\end{equation}
We see that the function $\phi$ is a solution of a
special kind of functional equation. Such equations are called
{\em refinement equations}, and their solutions are called  {\em refinable functions}
\footnote{Usually the terms ``refinable function'' and ``scaling function'' are
synonyms in the literature, and they are used in both senses: as a solution
to the refinable equation and as a function generating MRA.
We separate here the meanings of these terms.}.
It will be shown in Section~\ref{s3} that any test scaling function
(not necessary orthogonal) is  refinable.

A natural way for the construction of a MRA (see, e.g.,~\cite[\S 1.2]{NPS})
is the following. We start with a refinable function $\phi$
and define the spaces $V_j$ by~(\ref{17}).
It is clear that axioms (d) and (e) of Definition~\ref{de1} are fulfilled.
Of course, not any such function $\phi$ provides axiom $(a)$.
In the  real setting, the relation $V_0\subset V_{1}$ holds
if and only if the refinable function satisfies a refinement equation.
The situation is different in the $p$-adic case. Generally speaking, a refinement
equation (\ref{62.0-2*})  does not imply the including property
$V_0\subset V_{1}$ because the set of shifts $I_p$  does not
form a group.
Indeed, we need all the functions $\phi(\cdot-b)$,
$b\in I_p$, to belong to the space $V_1$, i.e., the identities
$\phi(x-b)=\sum_{a\in I_p}\alpha_{a,b}\phi(p^{-1}x-a)$ should be
fulfilled for all $b\in I_p$. Since $p^{-1}b+a$ is not in $I_p$ in general,
we  can not state that $\phi(\cdot-b)$ belongs to  $V_1$
for all $b\in I_p$.
Nevertheless, we will see below that a wide class of refinable equations
provide the including property.

Providing axiom (a) is a key moment for the construction of MRA.
Axioms (b) and (c) are fulfilled for a wide class of functions $\phi$
because of the following statements.

\begin{theorem}
\label{th1-2*}
If $\phi \in  L^2(\bQ_p)$ and $\widehat\phi$ is compactly supported,
then axiom $(c)$ of Definition~{\rm\ref{de1}} holds for the spaces
$V_j$ defined by~(\ref{17}).
\end{theorem}

\begin{proof} Let $\widehat\phi\subset B_M(0)$, $M\in\z$.
Assume that a function $f\in L^2(\bQ_p)$ belongs to any space $V_j$, $j\in\z$.
 Given  $j\in\n$ and $\epsilon>0$, there exists a  function
 $f_\epsilon(x):=\sum_{a\in I_p}\alpha_a\phi(p^jx-a)$, where the sum is finite,
  such that $\|f-f_\epsilon\|<\epsilon$. Using~(\ref{014}),
  it is not difficult to see that
  $\supp\widehat f_\epsilon\subset\supp\widehat\phi(p^{-j}\cdot)$,
  which yields that   $\widehat f_\epsilon(\xi)=0$ for any
  $\xi\not\in B_{M-j}(0)$.  Due to the Plancherel theorem,
  it follows that $\widehat f=0$ almost everywhere
  on $B_{M-j}(0)$. Since $j$ is an arbitrary positive integer,
 $\widehat f$ is equivalent to zero on $Q_p$.
\end{proof}

Another sufficient condition for axiom (c) was given in~\cite{Kh-Sh-S}:

\begin{theorem}
\label{th1-4*}
If $\phi \in  L^2(\bQ_p)$ and the system $\{\phi(x-a):a\in I_p\}$ is orthonormal,
then axiom $(c)$ of Definition~{\rm\ref{de1}} holds for the spaces
$V_j$ defined by~(\ref{17}).
\end{theorem}

\begin{theorem}
Let $\phi \in  L^2(\bQ_p)$, the spaces $V_j$,
$j\in \z$, be defined by~(\ref{17}),
and let $\phi(\cdot-b)\in\cup_{j\in \z} V_j$ for any $b\in Q_p$.
Axiom $(b)$ of Definition~{\rm\ref{de1}} holds for the spaces
$V_j$, $j\in\z$,   if and only if
\be
\bigcup\limits_{j\in\z}{\rm supp\,}\widehat\phi(p^{j}\cdot)=\bQ_p.
\label{dnn14}
\ee
\label{th1-3*}
\end{theorem}

\begin{remark}
It is not difficult to see that the assumption
$\phi(\cdot-b)\in\cup_{j\in \z} V_j$ for any $b\in Q_p$
is fulfilled whenever $\phi$ is a refinable function
and $\widehat\phi\subset B_0(0)$. We will see that
this assumption is also valid for a wide class of refinable
functions  $\phi$ for which $\widehat\phi\not\subset B_0(0)$.
\end{remark}

\begin{proof}
First of all we show that the space
$\overline {\cup_{j\in {\z}} V_j}$ is invariant with respect to all
shifts. Let $f\in \cup_{j\in {\z}} V_j$, $b\in\bQ_p$. Evidently,
$\phi(p^{-k}{\cdot}-t)\in\cup_{j\in \z} V_j$ for any $t\in Q_p$
and for any $k\in\z$.
Since the $L^2$-norm is invariant with
respect to the shifts,  it follows  that
$f(\cdot -b) \in \overline{\cup_{j\in\z} V_j}$.
 If now $g\in\overline{\cup_{j\in\z} V_j}$, then
approximating $g$ by the functions  $f\in \cup _{j\in\z} V_j$,
again using  the invariance of $L^2$-norm with respect to the shifts
, we derive $g(\cdot -b) \in \overline{\cup_{j\in\z} V_j}$.

For $X\subset L^2(\bQ_p)$, set $\widehat X=\{\ \widehat f: f\in X\}$.
By the Wiener theorem for $L^2$ (see, e.g., \cite{NPS}; all the arguments
of the proof given there may be repeated word for word with replacing
 $\bR$ by $\bQ_p$), a closed subspace
$X$ of the space $L^2(\bQ_p)$ is invariant with respect to the shifts
if and only if $\widehat X=L^2(\Omega)$ for some set
$\Omega\subset\bQ_p$. If now $X=\overline{\cup_{j\in\z}V_j}$, then
$\widehat X=L^2(\Omega)$. Thus $X=L^2(\bQ_p)$ if and only if $\Omega=\bQ_p$. Set
$\phi_j=\phi(p^{-j}\cdot),\ \ \Omega_0=\cup_{j\in\z}{\rm supp}\, \widehat\phi_j$
and prove that $\Omega=\Omega_0$. Since $\phi_j\in V_j,$ $j\in\z$,
we have ${\rm supp}\, \widehat\phi_j\subset\Omega$, and hence
$\Omega_0 \subset \Omega$.
Now assume that $\Omega\backslash\Omega_0$ contains a set of
positive measure $\Omega_1$. Let $f\in V_j$. Given   $\epsilon>0$,
there exists a  function
 $f_\epsilon(x):=\sum_{a\in I_p}\alpha_a\phi(p^jx-a)$, where the sum is finite,
  such that $\|f-f_\epsilon\|<\epsilon$. Using~(\ref{014}),
  we see that $\supp\widehat f_\epsilon\subset\supp\widehat\phi(p^{-j}\cdot)$,
  which yields that   $\widehat f_\epsilon(\xi)=0$ for any
  $\xi\not\in \Omega_1$.  Due to the Plancherel theorem,
  it follows that $\widehat f=0$ almost everywhere  on $\Omega_1.$
Hence the same is true for any $f\in \cup _{j\in\z} V_j$.
Passing to the limit we deduce that that the Fourier transform
of any $f\in X$ is equal to zero almost everywhere on  $\Omega_1$, i.e.,
$L^2(\Omega)=L^2(\Omega_0)$.
It remains to note that
${\rm supp\,} \widehat\phi_j={\rm supp\,}\widehat\phi({p}^{j}\cdot) $
\end{proof}

A real analog of Theorem~\ref{th1-3*} was proved by
 C.~de~Boor,  R.~DeVore and A.~Ron in~\cite5.

\section{Refinable functions}
\label{s3}

We are going to study $p$-adic refinable functions $\phi$. Let us restrict
ourselves to the consideration of  $\phi\in {\cD}$.
Evidently, each $\phi\in {\cD}$ is a
$p^M$-periodic  function for some $M\in \z$.
Denote by ${\cD}_N^M$ the set of all $p^M$-periodic  functions supported on $B_N(0)$.
Taking the Fourier transform
of the equality $\phi(x-p^M)=\phi(x)$, we obtain
$\chi_p(p^M\xi)\widehat\phi(\xi)=\widehat\phi(\xi)$, which holds for all $\xi$
if and only if $\supp\widehat\phi\subset B_M(0)$. Thus, the set ${\cD}_N^M$
consists of all locally constant functions $\phi$ such that $\supp\phi\subset B_N(0)$,
 $\supp\widehat\phi\subset B_M(0)$.

\begin{proposition}
Let $\phi,\psi \in  L^2(\bQ_p)$, $\supp \phi,\supp \psi\subset B_N(0)$, $N\ge0$,
 and let $b\in I_p$, $|b|_p\le p^N$. If
\be
\psi(\cdot-b)\in \overline{\span\{\phi(p^{-1}x-a),\ a\in I_p\}}
\label{19}
\ee
then
\be
\psi(x-b)=\sum_{k=0}^{p^{N+1}-1}h^\psi_{k,b}\phi\Big(\frac{x}{p}-\frac{k}{p^{N+1}}\Big)
\ \ \ \forall x\in Q_p.
\label{18}
\ee
\label{p1}
\end{proposition}

\begin{proof}
Given  $\epsilon>0$, there exist  functions
 $$
 f_\epsilon(x):=\sum_{a\in I_p\atop |a|_p\le p^{N+1}}\alpha_a\phi(p^jx-a),\ \ \
g_\epsilon(x):=\sum_{a\in I_p\atop |a|_p> p^{N+1}}\alpha_a\phi(p^jx-a),
$$
  where the sums are finite,
  such that $\|\psi(\cdot-b)-f_\epsilon-g_\epsilon\|<\epsilon$.
  If $x\in B_N(0)$, $|a|_p>p^{N+1}$, then
$|p^{-1}x-a|_p>p^{N+1}$ and hence $\phi(p^{-1}x-a)=0$. So, $g_\epsilon(x)=0$
whenever $x\in B_N(0)$. If $x\not\in B_N(0)$, then $\phi(x-b)=0$ and
$\phi(p^{-1}x-a)=0$ for all $a\in I_p$, $|a|_p\le p^{N+1}$.
So, $\phi(\cdot-b)-f_\epsilon(x)=0$ whenever $x\not\in B_N(0)$.
It follows that
$$
\|\psi(\cdot-b)-f_\epsilon\|^2=\int\limits_{B_N(0)}|\psi(x-b)-f_\epsilon|^2\,dx=
\int\limits_{B_N(0)}|\psi(x-b)-f_\epsilon-g_\epsilon|^2\,dx\le \epsilon^2.
$$
Hence
$$
\psi(\cdot-b)\in
\overline{\span\{\phi(p^{-1}x-a),\ a\in I_p,\ |a|_p\le p^{N+1}\}},
$$
which implies~(\ref{18}).
\end{proof}

\begin{corollary}
If $\phi \in  L^2(\bQ_p)$ is a refinable  function and
 $\supp \phi\subset B_N(0)$, $N\ge0$,  then  its refinement equation is
\begin{equation}
\label{62.0-5}
\phi(x)=\sum_{k=0}^{p^{N+1}-1}h_{k}\phi\Big(\frac{x}{p}-\frac{k}{p^{N+1}}\Big)
\ \ \ \forall x\in Q_p.
\end{equation}
\label{c3}
\end{corollary}

The proof immediately follows from Proposition~\ref{p1}.

\begin{corollary}
Let $\phi \in  L^2(\bQ_p)$ be a scaling  function of a MRA.
If  $\supp \phi\subset B_N(0)$, $N\ge0$,  then  $\phi$ is a refinable
function satisfying~(\ref{62.0-5}).
\label{c1}
\end{corollary}

The proof follows by combining axiom $(a)$ of Definition~{\rm\ref{de1}} with
Proposition~\ref{p1}.

Taking the Fourier transform of~(\ref{62.0-5})  and using
(\ref{014}),  we can rewrite the refinable equation in the form
\begin{equation}
\label{62.0-6}
{\widehat\phi}(\xi)=m_0\Big(\frac{\xi}{p^{N}}\Big){\widehat\phi}(p\xi),
\end{equation}
where
\begin{equation}
\label{62.0-7-1}
m_0(\xi)=\frac{1}{p}\sum_{k=0}^{p^{N+1}-1}h_{k}\chi_p(k\xi)
\end{equation}
is a trigonometric polynomial. It is clear that $m_0(0)=1$ whenever $\widehat\phi(0)\ne0$.

\begin{proposition}
\label{pr1-1}
If $\phi\in  L^2(\bQ_p)$ is a solution of refinable
equation~(\ref{62.0-5}), ${\widehat\phi}(0)\ne 0$, ${\widehat\phi}(\xi)$ is continuous at
the point $0$, then
\begin{equation}
\label{62.0-8}
{\widehat\phi}(\xi)={\widehat\phi}(0)\prod_{j=0}^{\infty}m_0\Big(\frac{\xi}{p^{N-j}}\Big).
\end{equation}
\end{proposition}

\begin{proof}
Since~(\ref{62.0-5}) implies~(\ref{62.0-6}), after
iterating {\rm(\ref{62.0-6})} $J$ times, $J\ge 1$, we have
$$
{\widehat\phi}(\xi)=\prod_{j=0}^{J}m_0\Big(\frac{\xi}{p^{N-j}}\Big)
{\widehat\phi}(p^{J}\xi).
$$
Taking into account that ${\widehat\phi}(\xi)$ is continuous at
the point $0$ and the fact that $|p^{N}\xi|_p=p^{-N}|\xi|_p\to 0$
as $N\to +\infty$ for any $\xi\in\bQ_p$, we obtain {\rm(\ref{62.0-8})}.
\end{proof}

\begin{corollary}
If $\phi \in  {\cD}_N^M$ is a refinable function, $N\ge0$, and ${\widehat\phi}(0)\ne 0$, then
(\ref{62.0-8}) holds.
\label{c2}
\end{corollary}

This statement follows immediately from Corollary~\ref{c1} and
Proposition~\ref{pr1-1}.

\begin{lemma}
Let
$
{\widehat\phi}(\xi)=C\prod_{j=0}^{\infty}
m_0\Big(\frac{\xi}{p^{N-j}}\Big),
$
where $m_0$ is a trigonometric polynomial with $m_0(0)=1$
and $C\in\r$.
If $\supp\widehat\phi\subset B_M(0)$,
then there exist at most
$\frac{\deg m_0}{p-1}$ integers $n$ such that
$0\le n<p^{M+N}$ and $\widehat\phi\lll\frac{n}{p^M}\rrr\ne0$.
\label{l1}
\end{lemma}

\begin{proof}
First of all we note that $\widehat\phi$ is a $p^N$-periodic
function satisfying~(\ref{62.0-6}).
Denote by $O_p$ the set of positive integers not divisible by~$p$.
Since $\supp\widehat\phi\subset B_M(0)$,
we have $\widehat\phi\lll\frac{k}{p^{M+1}}\rrr=0$ for all
$k\in O_p$
. By the definition of $\widehat\phi$ the equality
$\widehat\phi\lll\frac{k}{p^{M+1}}\rrr=0$ holds if and only if
there exists $\nu=1-N\ddd M+1$ such that $m_0\Big(\frac{k}{p^{N+\nu}}\Big)=0$.
Set
$$
\sigma_\nu:=\left\{l\in O_p:\ l<p^{N+\nu},
m_0\Big(\frac{l}{p^{N+\nu}}\Big)=0,\
m_0\Big(\frac{l}{p^{N+\mu}}\Big)\ne0\ \forall \mu=1-N\ddd\nu-1\right\},
$$
$v_\nu:=\sharp\,\sigma_\nu$. Evidently, $\sigma_\nu\subset O_p^{\prime}$
for all $\nu$, where $O_p^{\prime}=\{k\in O_p:\ k<p^{M+N+1}\}$,
and $\sigma_{\nu^{\prime}}\cap\sigma_{\nu}=\emptyset$ whenever
$\nu^{\prime}\ne\nu$. If $\widehat\phi\lll\frac{k}{p^{M+1}}\rrr=0$ for some
$k\in O_p$, then there exist a unique $\nu=1-N\ddd M+1$
and a unique $l\in\sigma_\nu$ such that $k\equiv l\pmod{p^{N+\nu}}$. Moreover,
for any $l\in\sigma_\nu$ there are exactly $p^{M-\nu+1}$ integers $k\in O_p^{\prime}$
(including~$l$) satisfying the above comparison.
It follows that
\be
\sum\limits_{\nu=1-N}^{M+1} p^{M-\nu+1}v_\nu=\sharp\,O_p^{\prime}=p^{M+N}(p-1).
\label{10}
\ee
Now if $l\in\sigma_\nu$, $\nu\le M$, then
$\widehat\phi\lll\frac{p^\gamma k}{p^{M}}\rrr=0$
 for all $\gamma=0,1\ddd M-\nu$, $k=l+rp^{N+\nu}$, $r=0,1\ddd p^{M-\nu-\gamma}-1$,
 i.e., each $l\in\sigma_\nu$ generates at least $1+p+\dots+p^{M-\nu}$
 distinct positive integers $n<p^{M+N}$ for which
 $\widehat\phi\lll\frac{n}{p^{M}}\rrr=0$. Hence
 \ban
v:= \sharp\,\left\{n:\ n=0,1\ddd p^{M+N}-1,
\widehat\phi\lll\frac{n}{p^{M}}\rrr=0\right\}\ge
\\
\sum\limits_{\nu=1-N}^{M}(1+p+\dots+p^{M-\nu})v_\nu=
\frac{1}{p-1}\sum\limits_{\nu=1-N}^{M}(p^{M-\nu+1}-1)v_\nu=
\\
\frac{1}{p-1}\sum\limits_{\nu=1-N}^{M+1}(p^{M-\nu+1}-1)v_\nu.
\ean
Since $\sum\limits_{\nu=1-N}^{M+1}v_\nu\le \deg m_0$, by using~(\ref{10}),
we obtain
$$
v\ge\frac{1}{p-1}\lll\sum\limits_{\nu=1-N}^{M+1}p^{M-\nu+1}v_\nu-
\deg m_0\rrr\ge p^{M+N}-\frac{\deg m_0}{p-1}.
$$
\end{proof}

For each $\phi\in{\cD}_N^M$, $M,N\ge0$, we assign the set
\be
L_\phi=\left\{l=0,1\ddd p^{M+N}-1:
\widehat\phi\lll\frac{l}{p^M}\rrr\ne0\right\}
\label{50}
\ee

\begin{theorem}
\label{t1}
Let $\phi\in{\cD}_N^M$, $M,N\ge0$ and ${\widehat\phi}(0)\ne 0$. If
\be
\phi(\cdot-b)\in \overline{\span\{\phi(p^{-1}x-a),\ a\in I_p\}}
\label{27}
\ee
for all $b\in I_p$, $|b|_p\le p^N$, then $\sharp\, L_\phi\le p^N$.
\end{theorem}

\begin{proof}
Let $b\in I_p$, $|b|_p\le p^N$. Because of Proposition~\ref{p1},
we can rewrite~(\ref{27}) in the form
$$
\phi(x-b)=\sum_{k=0}^{p^{N+1}-1}h_{k,b}\phi\Big(\frac{x}{p}-\frac{k}{p^{N+1}}\Big)
\ \ \ \forall x\in \bQ_p.
$$
Taking the Fourier transform, we obtain
\begin{equation}
\label{12}
{\widehat\phi}(\xi)\chi_p(b\xi)=m_b\Big(\frac{\xi}{p^{N}}\Big){\widehat\phi}(p\xi),
\ \ \ \forall \xi\in \bQ_p,
\end{equation}
where $m_b$ is a trigonometric polynomial, $\deg m_b<p^{N+1}$.
Combining~(\ref{12}) for $b=0$ with~(\ref{12}) for arbitrary $b$,
we obtain
$$
{\widehat\phi}(p\xi)\lll m_0\Big(\frac{\xi}{p^{N}}\Big)\chi_p(b\xi)-
 m_b\Big(\frac{\xi}{p^{N}}\Big)\rrr=0
\ \ \ \forall \xi\in \bQ_p,
$$
which is equivalent to
\be
F(\xi):={\widehat\phi}(p^{N+1}\xi)\lll m_0(\xi)\chi_p(p^{N}b\xi)-
 m_b(\xi)\rrr=0
\ \ \ \forall \xi\in \bQ_p.
\label{13}
\ee
Since  $\supp F\subset B_{M+N+1}(0)$ and $F$ is a $1$-periodic function,
(\ref{13}) holds if and only if $\widehat\phi\lll\frac{l}{p^{M+N+1}}\rrr=0$,
$l=0,1\ddd p^{M+N+1}-1$.

First suppose that $\deg m_0\ge p^N(p-1)$, i.e.,
$$
m_0(\xi)=\sum\limits_{k=0}^{K}h_k\chi_p(k\xi),\ \ \  h_K\ne0,
$$
where $K=K_Np^N+K_{N-1}p^{N-1}+\dots+K_0$,
$K_j\in D_p$, $j=0,1\ddd N$, $K_N=p-1$ (indeed,
if $K_N<p-1$, then $\deg m_0= K\le (p-2)p^N+(p-1)(1+p+\dots+p^{N-1})=
p^{N+1}-p^N-1<p^N(p-1)$). Set $b:= p-p^{-N}K$. It is not difficult to see that
$b\in I_p$, $|b|_p\le p^N$ and $K+bp^N=p^{N+1}$. We see that
the degree of the polynomial $t(\xi):=m_0(\xi)\chi_p(p^{N}b\xi)- m_b(\xi)$
is exactly $p^{N+1}$, and hence there exist at most $p^{N+1}$ integers $l$
such that $0\le l<p^{M+N+1}$, $t\lll\frac{l}{p^{M+N+1}}\rrr=0$. Thus,
$$
\sharp\,\left\{l:\ l=0,1\ddd p^{M+N+1}-1,
\widehat\phi\lll\frac{l}{p^{M}}\rrr=0\right\}\ge p^{M+N+1}-p^{N+1}.
$$
Taking into account that $\widehat\phi$ is a $p^N$-periodic function,
we obtain
\be
 \sharp\,\left\{l:\ l=0,1\ddd p^{M+N}-1,
\widehat\phi\lll\frac{l}{p^{M}}\rrr=0\right\}\ge p^{M+N}-p^{N}.
\label{14}
\ee
It remains to note that~(\ref{14}) is also fulfilled  whenever
$\deg m_0<p^N(p-1)$ because of Lemma~\ref{l1} and Corollary~\ref{c2}.
\end{proof}

\begin{theorem}
\label{t2} Let $\phi\in{\cD}_N^M$,  $M,N\ge0$,
$\sharp\,L_\phi\le p^N$, then
\be
\phi(x-b)=\sum_{a\in I_p}\alpha_{a,b}\phi(x-a) \ \ \
\forall b\in \bQ_p,
\label{11}
\ee
where the sum  is finite.
\end{theorem}

\begin{proof}
First we assume that $b\in Q_p$, $|b|_p\le p^{N}$,  and
prove that
\be
\phi(x-b)=\sum_{k=0}^{p^{N}-1}\alpha_{k,b}\phi\Big({x}-\frac{k}{p^{N}}\Big)
\ \ \ \forall x\in \bQ_p.
\label{29}
\ee
 Taking the
Fourier transform, we reduce~(\ref{29}) to
\begin{equation}
\label{12}
{\widehat\phi}(\xi)\chi_p(b\xi)=m_b\Big(\frac{\xi}{p^{N}}\Big){\widehat\phi}(\xi),
\ \ \ \forall \xi\in \bQ_p,
\end{equation}
where $m_b$ is a trigonometric polynomial, $\deg m_b<p^{N}$,
which is equivalent to
\be
f(\xi):={\widehat\phi}(p^{N}\xi)\lll
\chi_p(p^{N}b\xi)-
 m_b(\xi)\rrr=0
\quad \forall \xi\in \bQ_p.
\label{13}
\ee
Since  $\supp f\subset B_{M+N}(0)$
and $f$ is a $1$-periodic function, (\ref{13}) is equivalent to
$$
f\lll\frac{l}{p^{M+N}}\rrr=0, \forall l=0,1\ddd p^{M+N}-1,
$$
which holds if and only if
\ba
m_b\lll\frac{l}{p^{M+N}}\rrr=
\chi_p\lll\frac{bl}{p^{M}}\rrr,
\quad\forall l\in L_\phi.
\label{15}
\ea
  Hence we can find $m_b$ by
solving the linear system~(\ref{15}) with respect to the unknown
coefficients of $m_b$.  So, we proved~(\ref{12}), and
hence~(\ref{29}).

Next let  $b\in Q_p$, $|b|_p= p^{N+1}$, i.e.,
$b=\frac{b_{N+1}}{p^{N+1}}+b^\prime$, $b_{N+1}\in D_p$, $b_{N+1}\ne0$,
$|b^\prime|_p\le p^{N}$.
Using~(\ref{29}) with $b=b^\prime$, we have
$$
\phi(x-b)=\sum_{k=0}^{p^{N}-1}\alpha_{k,b^\prime}
\phi\Big({x}-\frac{k}{p^{N}}-\frac{b_{N+1}}{p^{N+1}}\Big)=
\sum_{k=0}^{p^{N}-1}\alpha_{k,b^\prime}
\phi\Big({x}-\frac{pk+b_{N+1}}{p^{N+1}}\Big).
$$
Taking into account that
$$
pk+b_{N+1}\le p(p^{N}-1)+(p-1)=p^{N+1}-1,
$$
we derive
$$
\phi(x-b)=\sum_{k=0}^{p^{N+1}-1}\alpha_{k,b}
\phi\Big({x}-\frac{k}{p^{N+1}}\Big)
\ \ \ \forall x\in \bQ_p.
$$
Similarly, we can prove by induction on $n$ that
$$
\phi(x-b)=\sum_{k=0}^{p^{N+n}-1}
\alpha_{k,b}\phi\Big(\frac{x}{p}-\frac{k}{p^{N+n}}\Big)
\ \ \ \forall x\in \bQ_p,
$$
whenever $b\in \bQ_p$, $|b|_p= p^{N+n}$.
\end{proof}

As a consequence we have the following statements.

\begin{corollary}
Let $\phi \in  {\cD}_N^M$ be a refinable function, $M,N\ge0$,
$L_\phi\le p^N$, and let  the spaces $V_j$ be defined by~(\ref{17}).
Then axiom (a) of Definition~\ref{de1}.
holds.
\label{c4}
\end{corollary}

\begin{corollary}
If a test function $\phi$ with ${\widehat\phi}(0)\ne 0$ generates a MRA,
then the corresponding spaces $V_j$, $j\in\z$, are invariant
with respect to all translations.
\label{c5}
\end{corollary}

\begin{theorem}
\label{t3}
A function $\phi\in{\cD}_N^M$, $M,N\ge0$,  with
${\widehat\phi}(0)\ne 0$ generates a MRA if and only if

(1) $\phi$ is refinable;

(2) 
there exist at most
$p^N$ integers $l$ such that $0\le l<p^{M+N}$ and
$\widehat\phi\lll\frac{l}{p^{M}}\rrr\ne0$.
\end{theorem}

\begin{proof}
If $\phi$ is a scaling function of a MRA, then (1) follows from
Corollary~\ref{c1}, and (2) follows from (1) and Theorem~\ref{t1}.

Now let conditions (1), (2) be fulfilled. Define the spaces $V_j$, $j\in\z$,
by~(\ref{17}). Axioms (d) and (e), evidently, hold.
Axiom (a) follows from Corollary~\ref{c4}. Axiom (b) follows from
Theorems~\ref{t2} and ~\ref{th1-3*}. Axiom (c) follows from
Theorems~\ref{th1-2*}.
\end{proof}

\begin{example}
Let $p=2$, $N=2$, $M=1$  $\phi$ be defined by (\ref{62.0-8}),
where $\widehat\phi(0)\ne0$, $m_0$ is given by (\ref{62.0-7-1}),
$m_0(1/4)=m_0(3/8)=m_0(7/16)=m_0(15/16)=0$ and $m_0(0)=1$.
It is not difficult to see that ${\rm supp\,}\widehat\phi\subset B_1(0)$,
${\rm supp\,}\widehat\phi\not\subset B_0(0)$ and
$\widehat\phi\Big(\frac{1}{2}\Big)=\widehat\phi\Big(\frac{3}{2}\Big)=
\widehat\phi\Big(\frac{5}{2}\Big)=
\widehat\phi(1)=0$, i.e, all the assumptions of Theorem~\ref{t3}
are fulfilled.
\label{e1}
\end{example}

\section{Orthogonal scaling functions}
\label{s4}

Now we are going to describe all orthogonal
scaling functions  $\phi\in{\cD}_N^M$.

\begin{theorem}
\label{th1-5*}
Let  $\phi\in{\cD}_N^M$,  $M,N\ge0$.
If    $\{\phi(x-a):a\in I_p\}$
is an orthonormal system, then
\be
\sum_{l=0}^{p^{M+N}-1}\left|{\widehat\phi}\lll\frac{l}{p^M}\rrr\right|^2
\chi_p\lll\frac{lk}{p^{M+N}}\rrr=p^N\delta_{k 0},
\quad k=0,{1},\dots,{p^{N}-1}.
\label{20}
\ee
\end{theorem}

\begin{proof}
Let $a\in I_p$. Due to the orthonormality of $\{\phi(x-a):a\in I_p\}$,
using the Plancherel theorem, we have
$$
\delta_{a 0}=
\langle\phi(\cdot),\phi(\cdot-a)\rangle
\int\limits_{\bQ_p}\phi(x)\overline{\phi(x-a)}\,dx
=\int\limits_{B_M(0)}|{\widehat\phi}(\xi)|^2\chi_p(a\xi)\,d\xi.
$$
Let $\xi\in B_M(0)$. There exists a unique $l=0,1\ddd p^{M+N}-1$
such that  $\xi\in B_{-N}\lll b_l\rrr$, $b_l=\frac{l}{p^M}$. It follows that
\ban
\int\limits_{B_M(0)}|{\widehat\phi}(\xi)|^2\chi_p(a\xi)\,d\xi
=\sum_{k=0}^{p^{M+N}-1}\int\limits_{|\xi-b_l|_p\le p^{-N}}
|{\widehat\phi}(\xi)|^2\chi_p(a\xi)\,d\xi\qquad\qquad\qquad\qquad
\\
=\sum\limits_{l=0}^{p^{M+N}-1}|{\widehat\phi}(b_l)|^2
\int\limits_{|\xi-b_l|_p\le p^{-N}}\chi_p(a\xi)\,d\xi
=\sum\limits_{l=0}^{p^{M+N}-1}|{\widehat\phi}(b_l)|^2\chi_p(ab_l)
\int\limits_{|\xi|_p\le p^{-N}}\chi_p(a\xi)\,d\xi
\\
=\frac{1}{p^{N}}\Omega(|p^{N}a|_p)
\sum_{l=0}^{p^{M+N}-1}|{\widehat\phi}(b_l)|^2\chi_p(ab_l).
\ean

To prove~(\ref{20}) it only remains to note that
 $\Omega(|p^{N}a|_p)= 0$ whenever $a\in I_p$,
 $p^Na\ne0,1\ddd p^N-1$.
\end{proof}

\begin{lemma}
Let $c_0,\ddd c_{n-1}$ be mutually distinct elements of the unit
circle $\{z\in\bC:\ |z|=1\}$. Suppose that there exist nonzero
reals $x_j$, $j=0,1\ddd n-1$, such that
\be
\sum_{j=0}^{n-1}c_j^kx_j=\delta_{k 0},\ \ \ k=0,1\ddd n-1.
\label{21}
\ee
Then $x_j=1/n$ for all $j$, and up to reordering
\be
c_j=c_0\ex{j/n},\ \ \ j=0,1\ddd n-1.
\label{22}
\ee
\label{l2}
\end{lemma}

\begin{proof}
In accordance with  Cramer's rule we have $x_j=\frac{\Delta_j}{\Delta}$,
$0\le j\le n-1$, where $\Delta=V(c)$ is the Vandermonde determinant
corresponding to $c=(c_0\ddd c_{N-1})$, and $\Delta_j$ is obtained from
$\Delta$ by replacing the $j$-th column with the transpose of the row
$(1,0\ddd 0)$. A straightforward computation shows that
$$
\Delta_j=(-1)^jV(c^{(j)})\prod\limits_{k\ne j}c_k,
$$
where $c^{(j)}$ is obtained from $c$ by removing the $j$-th coordinate. Thus,
\ba
x_j=(-1)^j\frac{V(c^{(j)})}{V(c)}\prod\limits_{k\ne j}c_k
=(-1)^j\prod\limits_{k\ne j}c_k\prod\limits_{k>l\atop k,l\ne j}
({c_k-c_l})\Big/{\prod\limits_{k>l}(c_k-c_l)}
\nonumber
\\
\prod\limits_{k\ne j}\frac{c_k}{c_k-c_j}=
\prod\limits_{k\ne j}\frac{1}{1-c^{-1}_kc_j}.
\label{35}
\ea
Next, for any $\alpha\in\r$, we have
\ban
1-e^{i\alpha}=2\sin\frac{\alpha}{2}
\lll\sin\frac{\alpha}{2}-i\cos\frac{\alpha}{2}\rrr=
2\sin\frac{\alpha}{2}e^{i\lll\frac{\alpha}{2}-\frac{\pi}{2}\rrr}.
\ean
Let us define $\alpha_j$, $j=0,1\ddd n-1$, by $c_j=e^{i\alpha_j}$.
Then from the above arguments and~(\ref{35}) it follows that
$$
x_j=\prod\limits_{k\ne j}\frac{1}{1-c^{-1}_kc_j}=
 e^{i\gamma}\sum\limits_{k\ne j}\lll2\sin\frac{\alpha_k-\alpha_j}{2}\rrr^{-1},
$$
where
$$
\gamma=\sum\limits_{k\ne j}\frac{\alpha_k-\alpha_j+\pi}{2}=
\theta-\frac{n}{2}\alpha_j,\ \ \ \theta=
\frac12\lll(n-1)\pi+\sum\limits_{k=0}^{n-1}\alpha_k\rrr
$$
By the lemma's hypothesis $x_j\in\r$, whence
$\gamma\equiv0\pmod{\pi}$ and consequently
$n\alpha_j\equiv2\theta\pmod{2\pi}$. Thus up to
reordering $\alpha_j=\alpha_0+\frac{2\pi j}{n}$,
which implies~(\ref{22}), and consequently that $x_j=1/n$ for all $j$.
\end{proof}

\begin{theorem}
\label{t4}
Let  $\phi\in{\cD}_N^M$ be an orthogonal
scaling function and  $\widehat\phi(0)\ne0$.
Then ${\rm supp\,{\widehat\phi}}\subset B_0(0)$.
\end{theorem}

\begin{proof}
Without loss of generality, we can assume that  $M,N\ge0$.
Combining Theorems~\ref{t3} and \ref{th1-5*}, we have
$$
\sum_{j=0}^{p^{N}-1}\left|{\widehat\phi}\lll\frac{l_j}{p^M}\rrr\right|^2
\chi_p\lll\frac{l_jk}{p^{M+N}}\rrr=p^N\delta_{k 0},
\quad k=0,{1},\dots,{p^{N}-1}.
$$
By Lemma~\ref{l2}, $l_j=l_0+jp^M$ and ${\widehat\phi}\lll\frac{l_j}{p^M}\rrr=1$.
Taking into account that  $\widehat\phi(0)\ne0$, we deduce
$l_0=0$, i.e., ${\widehat\phi}(j)=1$, $j=0,1\ddd p^{N}-1$.
Since $\widehat\phi$ is a $p^N$-periodic function,
it follows from Theorem~\ref{t3}  that
$\widehat\phi\lll\frac{l}{p^M}\rrr=0$ for all $l\in\z$
not divisible by $p^M$.
This yields ${\rm supp\,{\widehat\phi}}\subset B_0(0)$.
\end{proof}

So any test function $\phi$ generating a MRA belongs to the class
${\cD}_N^0$. All such functions were described in~\cite{Kh-Sh-S}.
The following theorem summarizes these results.

\begin{theorem}
\label{th1-6*}
Let ${\widehat\phi}$ be defined by {\rm(\ref{62.0-8})}, where $m_0$ is
the trigonometric polynomial~{\rm(\ref{62.0-7-1})} with $m_0(0)=1$.
If $m_0\big(\frac{k}{p^{N+1}}\big)=0$ for all $k=1,\dots,p^{N+1}-1$
not divisible by $p$, then  $\phi\in{\cD}_N^0$.
If, furthermore,
$\big|m_0\big(\frac{k}{p^{{N+1}}}\big)\big|=1$ for all $k=1,\dots,p^{N+1}-1$
divisible by $p$, then $\{\phi(x-a):a\in I_p\}$ is an orthonormal system.
Conversely, if  ${\rm supp\,{\widehat\phi}}\subset B_0(0)$
and the system $\{\phi(x-a):a\in I_p\}$
is orthonormal, then $\big|m_0\big(\frac{k}{p^{{N+1}}}\big)\big|=0$
whenever $k$ is not divisible
by $p$, $\big|m_0\big(\frac{k}{p^{{N+1}}}\big)\big|=1$ whenever $k$
is divisible by $p$, \ $k=1,2,\dots,p^{N+1}-1$,  and $|{\widehat\phi}(x)|=1$
for any $x\in B_0(0)$.
\end{theorem}

\begin{theorem}
There exists a unique MRA generated by an orthogonal scaling test function.
\label{t6}
\end{theorem}

\begin{proof}
Let  a MRA $\{V_j\}_{j\in\z}$  is generated by a
test scaling function $\phi$ such that the system
 $\{\phi(x-a):a\in I_p\}$ is  orthonormal. We prove that
 this MRA  coincides with the Haar MRA $\{V^H_j\}_{j\in\z}$
generated by the scaling function $\phi^H=\Omega(|\cdot|_p)$.
Evidently, it suffices to check that $V_0=V_0^H$.
Let $f\in V_0$. It follows from Theorem~\ref{t4} that
${\rm supp}\,{\widehat\phi}\subset B_0(0)$. Hence
${\rm supp}\,{\widehat f}\subset B_0(0)$, i.e.
$\widehat f\in L^2(B_0(0))$.
It is well known that each contineous character of the additive group of
the ring $\z_p$ is of the form $\chi_p(a\xi)$, ${a\in I_p}$.
Since this group is compact, by the Peter-Weyl theorem the
set of all these characters is an orthonormal basis for $L^2(B_N(0))=L^2(\z_p)$
(see, e.g.\cite{P}).Thus
we have $\widehat f(\xi)=\phi^H(\xi)\sum_{a\in I_p}\alpha_a\chi_p(a\xi)$,
$\sum_{a\in I_p}|\alpha_a|^2<\infty$. Taking the Fourier transform
and using~(\ref{014}) and (\ref{14.1}), we obtain
$f(x)=\sum_{a\in I_p}\phi^H(x-a)$. So $V_0\subset V_0^H$.
To prove the inclusion $V^H_0\subset V_0$ we will check that $\phi^H(\cdot-b)\in V_0$
for any $b\in I_p$. By Theorem~\ref{th1-6*}, the function $(\widehat\phi)^{-1}$
is bounded on $B_0(0)$. This yields that $(\widehat\phi)^{-1}\in L^2(B_N(0))$.
Hence,  $\chi_p(b\xi)(\widehat\phi(\xi))^{-1}=\phi^H(\xi)\sum_{a\in I_p}\beta_a\chi_p(a\xi)$,
$\sum_{a\in I_p}|\beta_a|^2<\infty$, which may be rewritten as
 $\chi_p(b\xi)\phi^H(\xi)=\widehat\phi(\xi)\sum_{a\in I_p}\beta_a\chi_p(a\xi)$.
Taking the Fourier transform
and using again~(\ref{014}), (\ref{14.1}), we obtain
$\phi^H(x-b)=\sum_{a\in I_p}\beta_a\phi(x-a)$.
\end{proof}

\section{Construction of wavelet frames}
\label{s4}

\begin{definition}
Let $H$ be a Hilbert space. A system $\{f_n\}_{n=1}^\infty\subset H$
is said to be a frame if there exist positive constants $A,B$
({\em frame boundaries}) such that
$$
A\|f\|^2\le \sum_{n=1}^\infty|\langle f,f_n\rangle|^2\le B\|f\|^2
\ \ \ \forall f\in H.
$$
\end{definition}

We are interested in the construction of $p$-adic wavelet frames, i.e.,
frames in $L^2(\bQ_p)$ consisting of functions $p^{j/2}\psi^{(\nu)}(p^{-j}x-a)$,
$a\in I_p$, $\nu=1\ddd r$.

Our general scheme of construction looks as follows. Let $\{V_j\}_{j\in\z}$
be a MRA. As above, we define
the wavelet space $W_j$, $j\in \z$,  as the orthogonal complement
of $V_j$ in $V_{j+1}$, i.e., $V_{j+1}=V_j\oplus W_j$.
 It is not difficult to see that
$f\in W_j$ if and only if $f(p^{j}\cdot)\in W_0$, and $W_j\perp
W_k$ whenever $j\ne k$. If now there exist functions
$\psi^{(\nu)}\in L^2(Q_p)$, $\nu=1\ddd r$, ({\em a set of wavelet
functions}) such that \be W_0=\overline{\span\{\psi^{(\nu)}(x-a),\
\nu=1\ddd r,\ a\in I_p\}}, \label{23} \ee then we have a wavelet
system \be \{p^{j/2}\psi^{(\nu)}(p^{-j}x-a),\ \nu=1\ddd r, a\in
I_p, j\in\z \}. \label{36} \ee It will be proved that such a
system is a frame in $L^2(Q_p)$ whenever  $\psi^{(\nu)}$ are
compactly supported functions.

\begin{theorem} Let  $\psi^{(\nu)}$, $\nu=1\ddd r$, be a set of
 compactly supported wavelet functions for a MRA  $\{V_j\}_{j\in\z}$. Then
the   system~(\ref{36})  is a frame in $L^2(Q_p)$.
\label{t5}
 \end{theorem}

\begin{proof} First we will prove that  the system
$\{\psi^{(\nu)}(\cdot-a),\ \nu=1\ddd r,\ a\in I_p\}$ is a frame in  the wavelet
space $W_0$. Let  $\supp\psi^{(\nu)}\subset B_N(0)$, $\nu=1\ddd r$, $N\ge0$. Set
Ïîëîæèì $a^{n, l}=\frac{l}{p^{N+n}}$, $l\in L(n)$, where  $L(n)$ is the set
of integers $l$, $0\le l< p^n$, which are not  divisible by $p$,
\ban
W_0^0&=&{\span\{\psi^{(\nu)}(x-a):\ \nu=1\ddd r,\  a\in  I_p\cap B_N(0)\}},
\\
W_0^{n, l}&=&{\span\{\psi^{(\nu)}(x-a):\ \nu=1\ddd r,\  a\in  I_p\cap  B_N(a^{n, l})\}},
\quad n\in\n,\quad  l\in L(n).
\ean
Since the disks $B_N(0)$, $B_N(a^{n, l})$  are mutually disjoint and the
union  of them is $\bQ_p$, each function $f\in W_0$ may be represented in the form
 $$
 f=f^0+\sum\limits_{n=1}^\infty\sum_{l\in L(n)}f^{n,l},
 \ \ f^0=f\Big|_{B_N(0)}, f^{n,l}=f\Big|_{B_{N}(a^{n,l})}.
 $$
Due to~(\ref{23}),  given $\epsilon>0$, there exists a sum
 $\sum\limits_{a\in I_p}\sum\limits_{\nu=1}^{r}\alpha_a\psi^{(\nu)}(x-a)=:f_\epsilon(x)$,
 such that
$\|f-f_\epsilon\|<\epsilon$.  If $x\in  B_N(0)$, then
$f_\epsilon(x)=
\sum\limits_{a\in I_p\cap  B_N(0)}
\sum\limits_{\nu=1}^{r}\alpha_a\psi^{(\nu)}(x-a)=:f^0_\epsilon(x)$.
Since $\supp f^0\subset B_N(0)$, $\supp f_\epsilon^0\subset B_N(0)$, we have
$$
\|f-f_\epsilon\|^2\ge\int\limits_{B_N(0)}|f-f_\epsilon|^2
=\int\limits_{B_N(0)}|f^0-f^0_\epsilon|^2=\|f^0-f^0_\epsilon\|^2.
$$
Hence, $f^0\in W_0^0$. Similarly,
$f^{n, l}\in W_0^{n, l}$. It is not difficult to see
that the spaces $W_0^0$, $W_0^{n,l}$  are mutually orthogonal.
Thus we proved that
\be
W_0=W_0^0\oplus\lll \bigoplus\limits_{n=1}^\infty
\bigoplus\limits_{l\in L(n)}W_0^{n,l}\rrr.
\label{26}
\ee

Since $W_0^0$ is a finite
dimensional space and
$\{\psi^{(\nu)}(\cdot-a),\ \nu=1\ddd r,\  a\in  I_p\cap B_N(0)\}$
is a representing system for $W_0^0$, this system
is a frame. Hence there exist
 positive constants $A,B$ such that
\be
 A\|f^0\|^2\le\sum_{a\in I_p\cap B_N(0)}\sum\limits_{\nu=1}^{r}
 |\langle f^0,\psi^{(\nu)}(\cdot-a)\rangle|^2\le B\|f^0\|^2
  \ \ \ \forall f^0\in W_0^0.
 \label{25}
\ee
 If  $f^{n,l}\in W^{n,l}_0$, we have
\ban
\sum_{a\in I_p\cap B_{N}(a^{n,l})}\sum\limits_{\nu=1}^{r}|\langle
f^{n,l},\psi^{(\nu)}(\cdot-a)\rangle|^2= \sum_{a\in I_p\cap B_{N}(0)}
\sum\limits_{\nu=1}^{r}|\langle f^{n,l},\psi^{(\nu)}(\cdot-a^{n,l}-a)\rangle|^2=
\\
\sum_{a\in I_p\cap B_{N}(0)}\sum\limits_{\nu=1}^{r}|\langle
f^{n,l}(\cdot+a^{n,l}),\psi^{(\nu)}(\cdot-a)\rangle|^2.
\ean
Since $f^{n,l}(\cdot+a^{n,l})\in W_0^0$, it follows from~(\ref{25}) that
$$
 A\|f^{n,l}\|^2\le\sum_{a\in I_p\cap B_N(n,l)}\sum\limits_{\nu=1}^{r}
 |\langle f^{n,l},\psi^{(\nu)}(\cdot-a)\rangle|^2\le B\|f^{n,l}\|^2
  \ \ \ \forall f^{n,l}\in W_0^{n,l}.
$$
Taking into account~(\ref{26}), we  derive
\ban
  A\|f\|^2\le\sum_{a\in I_p}\sum\limits_{\nu=1}^{r}
  |\langle f,\psi^{(\nu)}(\cdot-a)\rangle|^2\le B\|f\|^2
 \ \ \ \forall f\in W_0.
\label{24}
\ean
So, we proved that the system $\{\psi^{(\nu)}(x-a),\ \nu=1\ddd r,\ a\in I_p\}$
is a frame in   $W_0$. Evidently, the system
$\{p^{j/2}\psi^{(\nu)}(p^{-j}x-a),\ \nu=1\ddd r, a\in I_p, \}$ is a frame in
$W_j$ with the same frame boundaries for any $j\in\z$. Since
${\bigoplus\limits_{j\in\z}W_j}= L^2(\bQ_p)$, it follows that the
union of these frames is a frame in $L^2(\bQ_p)$.
\end{proof}

Now we discuss how to construct a desirable set of wavelet
functions $\psi^{(\nu)}$, $\nu=1\ddd r$. Let a MRA
$\{V_j\}_{j\in\z}$  is generated by a scaling function
$\phi\in{\cD}_N^M$, $\widehat\phi(0)\ne0$. First of all we should
provide $\psi^{(\nu)}\in V_1$. Let us look for $\psi^{(\nu)}$ in
the form $$
\psi^{(\nu)}(x)=\sum_{k=0}^{p^{N+1}-1}g^{(\nu)}_{k}\phi\Big(\frac{x}{p}-\frac{k}{p^{N+1}}\Big)
$$ Taking the Fourier transform and using (\ref{014}),  we have $$
{\widehat\psi^{(\nu)}}(\xi)=n^{(\nu)}_0\Big(\frac{\xi}{p^{N}}\Big){\widehat\phi}(p\xi),\
\ \ $$ where $n^{(\nu)}_0$ is a trigonometric polynomial ({\em
wavelet mask}) given by $$
n^{(\nu)}_0(\xi)=\frac{1}{p}\sum_{k=0}^{p^{N+1}-1}g^{(\nu)}_{k}\chi_p(k\xi)
$$ Evidently, $\psi^{(\nu)}\in{\cD}_N^{M+1}$. By Theorem~\ref{t1},
there exist at least $p^{M+N}-p^N$ integers $l$ such that $0\le
l<p^{M+N}$, $\widehat\phi\lll\frac{l}{p^{M}}\rrr=0$. Choose
$n^{(\nu)}_0$ satisfying the following property: if $l\in L_\phi$,
i.e. $\widehat\phi\lll\frac{l}{p^{M}}\rrr\ne0$ for some $l=0,1\ddd
p^{M+N}-1$, then $n^{(\nu)}_0\lll\frac{l}{p^{M+N}}\rrr=0$. This
yields that $\widehat\psi^{(\nu)}\lll\frac{l}{p^{M}}\rrr=0$
whenever
 $0\le l<p^{M+N}$, $\widehat\phi\lll\frac{l}{p^{M}}\rrr\ne0$.

Let $a, b\in I_p$. Using the Plancherel theorem and the arguments
of Theorem~\ref{th1-5*}, we have
\ban
\langle\phi(\cdot-a),\psi^{(\nu)}(\cdot-b)\rangle=
\int\limits_{\bQ_p}\phi(x-a)\overline{\psi^{(\nu)}(x-b)}\,dx=
\\
\int\limits_{B_M(0)}{\widehat\phi}(\xi)
\overline{\widehat\psi^{(\nu)}(\xi)}\chi_p((b-a)\xi)\,d\xi=
\\
\sum_{l=0}^{p^{M+N}-1}\int\limits_{|\xi-p^{-M}l|_p\le p^{-N}}
{\widehat\phi}(\xi)
\overline{\widehat\psi^{(\nu)}(\xi)}\chi_p((b-a)\xi)\,d\xi=
\\
\sum\limits_{l=0}^{p^{M+N}-1}{\widehat\phi}\lll\frac{l}{p^M}\rrr
\overline{\widehat\psi^{(\nu)}\lll\frac{l}{p^M}\rrr}
\int\limits_{|\xi-p^{-M}l|_p\le p^{-N}}\chi_p(a\xi)\,d\xi=0.
\ean
It follows that $\overline{\span\{\psi^{(\nu)}(x-a),\ \nu=1\ddd r,\ a\in
I_p\}}\perp V_0$. On the other hand, due to Theorem~\ref{t2}, we
have $\overline{\span\{\psi^{(\nu)}(x-a),\ \nu=1\ddd r,\ a\in I_p\}}\subset
V_1$. Hence,
\be
\overline{\span\{\psi^{(\nu)}(x-a),\ \nu=1\ddd r,\ a\in
I_p\}}\subset W_0.
\label{30}
\ee
 It is clear from the proof of
Theorem~\ref{t2} that
\ba
&&\phi\lll x-\frac{l}{p^{N}}\rrr=
\sum_{k=0}^{p^{N+1}-1}h_{kl}\phi\Big(\frac{x}{p}-\frac{k}{p^{N+1}}\Big),
\ \ \ l=0\ddd p^N-1,
\label{32}
\\
&&\psi^{(\nu)}\lll
x-\frac{l}{p^{N}}\rrr=
\sum_{k=0}^{p^{N+1}-1}g^{(\nu)}_{kl}\phi\Big(\frac{x}{p}-\frac{k}{p^{N+1}}\Big),
\  l=0\ddd p^N-1,\ \nu=1\ddd r.
\label{31}
\ea
If the functions in the right hand side can be expressed as linear combinations
of the functions in the left hand side of ~(\ref{32}),~(\ref{31}),
i.e.
\ba
\span\left\{\phi\Big(\frac{x}{p}-a\Big), \ a\in I_p\cap B_{N+1}(0),
\right\} \subset
\label{51}
\\
{\span\{\phi(x-a),\psi^{(\nu)}(x-a),\ \nu=1\ddd r,\ a\in I_p\cap
B_{N}(0)\}}, \nonumber \ea then
$W_0\subset\overline{\span\{\psi^{(\nu)}(x-a),\ \nu=1\ddd r,\ a\in
I_p\}}.$ Taking into account~(\ref{30}), we deduce that
$\psi^{(\nu)},\ \nu=1\ddd r,$ is a set of wavelet functions.

Inclusion~(\ref{51}) is, evidently, fulfilled whenever the linear system
\ban
\sum_{k=0}^{p^{N+1}-1}h_{kl}x_k=0,
\ \ \ l=0\ddd p^N-1,
\\
\sum_{k=0}^{p^{N+1}-1}g^{(\nu)}_{kl}x_k,
\  l=0\ddd p^N-1,\ \nu=1\ddd r
\ean
has no non-trivial solutions. In particular, in the case $r=p-1$,
the system has no non-trivial solutions if and only if
the determinant is not equal zero.
It is not quite clear how to construct  functions
 $\psi^{(\nu)}$ providing~(\ref{51}) for arbitrary $\phi$ , but we
will show how to succeed in the case $\deg m_0\le(p-1) p^N$.
Such a masks with $p=2$ was presented in Example~\ref{e1}.

Assume that $\deg m_0\le(p-1) p^N$. In this case
$$
\phi\lll x\rrr=
\sum_{k=0}^{(p-1)p^{N}}h_{k}\phi\Big(\frac{x}{p}-\frac{k}{p^{N+1}}\Big).
$$

Define the wavelet masks
$n_0^{(\nu)}$, $\nu=1\ddd p-1$, by
\ban
n_0^{(\nu)}(\xi)=
\chi_p\lll(\nu-1)p^N\xi\rrr
(\chi_p(\xi)-1)^{p^N-\sharp\,L_\phi}\prod\limits_{l\in L_\phi}
\lll\chi_p\lll\xi\rrr-\chi_p\lll\frac{l}{p^{M+N}}\rrr\rrr=
\\
\frac{1}{p}\sum_{k=(\nu-1)p^N}^{\nu p^{N}}g_{k-(\nu-1)p^N}\chi_p(k\xi),
\quad \nu=1\ddd p-1,
\ean
(recall that $\sharp\, L_\phi\le p^N$ because of  Theorem~\ref{t1}).
So, system~(\ref{32}),~(\ref{31}) looks as
follows:
\ban
\phi\lll x-\frac{l}{p^{N}}\rrr=
\sum_{k=l}^{(p-1)p^{N}+l}h_{k-l}\phi\Big(\frac{x}{p}-\frac{k}{p^{N+1}}\Big)
\ \ \ l=0\ddd p^N-1,
\\
\psi^{(\nu)}\lll x-\frac{l}{p^{N}}\rrr=
\sum_{k=(\nu-1)p^N+l}^{\nu p^{N}+l}g_{k-(\nu-1)p^N-l}
\phi\Big(\frac{x}{p}-\frac{k}{p^{N+1}}\Big),\hspace{1cm}
\\
\nu=1\ddd p-1,\quad l=0\ddd p^N-1.
\ean
The determinant of the system  equals to
$$
\left|
\begin{array}{llllllllll}
 h_0 & h_1&\hdots& h_{p^N-1}& h_{p^N}&\dots&h_{(\nu-1)p^N}&0&\hdots&0
 \\
 0& h_0 &\dots& h_{p^{N}-2}& h_{p^{N}-1}&\hdots& h_{(\nu-1)p^N-1}&h_{(\nu-1)p^N}&\hdots&0
 \\
 \hdots&\hdots&\hdots&\hdots&\hdots&\hdots&\hdots&\hdots&\hdots&\hdots
 \\
  0& 0 &\dots& h_0& h_1& \dots&h_{(\nu-2)p^N+1}&h_{(\nu-2)p^N+2}&\hdots&h_{(\nu-1)p^N}
  \\
  g_0 & g_1&\hdots& g_{p^N-1}& g_{p^N}&\hdots&0&0&\hdots&0
 \\
 0& g_0 &\hdots& g_{p^{N}-2}& g_{p^{N}-1}&\hdots&0&0&\hdots&0
 \\
 \hdots&\hdots&\hdots&\hdots&\hdots&\hdots&\hdots&\hdots&\hdots&\hdots
 \\
  0& 0 &\hdots&0&0&\hdots& g_1&  g_2&\hdots&g_{p^N}
\end{array}
\right|
$$
This determinant  is so called resultant.
The resultant is not equal to zero if and only if the algebraic
polynomials with the coefficients $g_0,g_1\ddd g_{p^N}$ and
$h_0,h_1\ddd h_{p^N}$ respectively do not have joint zeros
 (see, e.g., \cite{L}).   But this holds
 because the trigonometric polynomials
 $m_0$ and $n^{(1)}_0$ do not have joint zeros by
 construction.

\bibliographystyle{amsplain}

\end{document}